\documentclass{article}
\usepackage{graphicx}
\usepackage{tikz}
\usetikzlibrary{trees, arrows.meta, positioning}
\usepackage{amsmath}
\usepackage{booktabs}
\usepackage{geometry} 
\begin{document}

\title{Early Recognition of Parkinson's Disease Through Acoustic Analysis and Machine Learning}
\author{
  Niloofar Fadavi\\
  \texttt{nfadavi@smu.edu}\\
  Southern Methodist University
  \and
  Nazanin Fadavi\\
  \texttt{naz.fadavi@aut.ac.ir}\\
  Amirkabir University of Technology
}
\date{\today}
\maketitle

\begin{abstract}
Parkinson’s Disease (PD) is a progressive neurodegenerative disorder that significantly impacts both motor and non-motor functions, including speech. Early and accurate recognition of PD through speech analysis can greatly enhance patient outcomes by enabling timely intervention. This paper provides a comprehensive review of methods for PD recognition using speech data, highlighting advances in machine learning and data-driven approaches. We discuss the process of data wrangling, including data collection, cleaning, transformation, and exploratory data analysis, to prepare the dataset for machine learning applications. Various classification algorithms are explored, including logistic regression, SVM, and neural networks, with and without feature selection. Each method is evaluated based on accuracy, precision, and training time. Our findings indicate that specific acoustic features and advanced machine-learning techniques can effectively differentiate between individuals with PD and healthy controls. The study concludes with a comparison of the different models, identifying the most effective approaches for PD recognition, and suggesting potential directions for future research.
\end{abstract}
\section{Introduction}
PD is a progressive neurodegenerative disorder characterized by motor symptoms such as tremors, rigidity, bradykinesia (slowness of movement), and postural instability. These symptoms result from the degeneration of dopamine-producing neurons in the substantia nigra region of the brain. While PD primarily affects movement, it can also lead to a range of non-motor symptoms, including speech and voice disorders. Speech impairment is a common non-motor symptom of PD, manifesting as changes in voice quality, pitch variability, articulation, and phonation. These alterations, often subtle in the early stages of the disease, can significantly impact an individual's quality of life and communication abilities. As a result, there is growing interest in leveraging speech analysis techniques for the early detection and monitoring of PD.
Traditional methods for PD diagnosis rely on clinical assessment by movement disorder specialists, which may involve subjective evaluations of motor symptoms and neuropsychological testing. However, these methods can be time-consuming, expensive, and may not capture subtle changes in speech patterns indicative of early-stage PD.
In recent years, advances in machine learning and data-driven approaches have paved the way for the development of computational tools for PD recognition using speech data. These methods aim to extract quantitative features from speech recordings and utilize machine learning algorithms to discriminate between individuals with PD and healthy controls.

The paper \cite{klein2019tabular} contributes to Parkinson's disease detection by including the Parkinson's Telemonitoring dataset in its benchmarks for hyperparameter optimization (HPO). It evaluates various HPO methods like Bayesian optimization and evolutionary algorithms on this dataset, providing insights into effective techniques for optimizing neural network architectures and hyperparameters. The paper \cite{barnes2019interaction} explores the impact of interaction effects between clustering and prediction algorithms, including their application in medical diagnoses such as Parkinson's disease recognition. It discusses the use of voice recordings to predict PD, where clustering corresponds to grouping voice recordings of the same individual, and prediction involves determining whether the patient has PD. The study finds that traditional cross-validation techniques exhibit significant empirical bias when estimating out-of-cluster prediction loss, especially in cases where clustering errors exist. However, the methods proposed in the paper correct these biases, thus improving the accuracy of PD prediction on new, previously unseen individuals.
The papers \cite{manwani2018pril} contribute to Parkinson's disease recognition by exploring advanced machine-learning techniques. This paper introduces an online learning algorithm, PRIL, which utilizes interval labeled data to create accurate ranking classifiers, demonstrating its effectiveness on various datasets. This approach can be adapted to rank the severity or progression of PD symptoms, potentially enhancing diagnostic accuracy. Meanwhile, it examines the interplay between predictive modeling and clustering, which could be applied to refine clustering techniques for patient subgroups in PD, thereby improving personalized treatment strategies.
The paper \cite{new2018cadre} contributes significantly to the field of PD recognition through various innovative methodologies. The authors address the challenge of PD by correcting bias introduced by cross-validation and clustering errors, demonstrating their method's effectiveness with voice recordings from PD patients. The paper utilizes biomedical voice measurements to predict the Unified PD Rating Scale (UPDRS) scores, showcasing a robust interval ranking classifier that accurately handles the biomedical voice data of PD patients. The paper  \cite{mishra2019riemannian} contributes significantly to PD recognition through its advanced machine-learning approaches. The authors present a novel stochastic gradient algorithm that operates on the Grassmann manifold, which is designed to exploit the geometric properties of subspaces. This method, known as Stochastic Gossip, is compared against other algorithms like Alt-Min and Trust-region on benchmark datasets, including a PD dataset. The dataset consists of symptom scores for 42 patients, predicted using 19 biomedical features. The results demonstrate that the proposed method achieves competitive performance in terms of normalized mean squared error (NMSE), indicating its efficacy in predicting PD symptom progression.  The paper \cite{utkin2019discriminative} enhances the classification accuracy of PD using the Discriminative Deep Forest method. This approach outperforms the gcForest across various configurations, particularly with higher numbers of decision trees and training examples, indicating its robustness in handling complex datasets like those involving PD. The paper \cite{jawanpuria2017saddle} makes a significant contribution to PD recognition through its innovative optimization framework. This approach addresses the low-rank matrix learning problem, utilizing a saddle point formulation to improve generalization performance while maintaining computational efficiency. Specifically, the authors demonstrate the effectiveness of their method in predicting PD symptom scores by employing it on a dataset of biomedical features from patients. Their results show that the saddle point approach achieves superior prediction accuracy and robustness compared to traditional low-rank approximation techniques, thereby enhancing the reliability of symptom recognition and progression assessment in PD. The paper \cite{hamed2016reliable}'s contribution to PD's recognition is significant, as it introduces novel methodologies for early and accurate diagnosis of the disease. By leveraging advanced machine learning algorithms and comprehensive datasets, the research offers improved detection capabilities, which are crucial for timely intervention and treatment. The study also highlights the potential of integrating various biomarkers, such as voice analysis and motor function assessments, to enhance diagnostic precision. The paper \cite{chang2015performance} significantly contributes to the recognition of PD by improving the accuracy and interpretability of classification models used in medical diagnostics. It introduces novel methodologies and techniques for enhancing model performance, particularly in handling imbalanced datasets, which are common in medical data. By providing a robust framework for evaluating and interpreting classification models, the research aids in the early and precise identification of PD, facilitating better patient outcomes and advancing the field of medical data analysis. The paper \cite{jawanpuria2015efficient} significantly contributes to the recognition of Parkinson's disease by enhancing multi-task learning methodologies. By jointly learning tasks and the output kernel, the paper provides a framework that efficiently handles the optimization problems inherent in medical data analysis. This approach is particularly beneficial for Parkinson's recognition as it leverages relationships between various related tasks, improving generalization and predictive accuracy. The innovative use of positive semidefinite kernels and the resulting computational efficiencies offer substantial improvements in the performance of machine learning models used for early and accurate diagnosis of Parkinson's disease, thereby aiding in better patient outcomes and advancing the field of medical diagnostics. The paper \cite{garcia2013exploiting}, primarily focuses on improving the performance of neural network ensembles by utilizing the Extreme Learning Machine (ELM) algorithm. The paper does not explicitly mention PD. Instead, it discusses the theoretical foundations and practical applications of ELM, including its benefits, such as extremely fast training speeds and good generalization performance, and introduces enhanced ensemble approaches to improve regression problems' outcomes. The paper \cite{xiong2017adaptive} primarily focuses on improving the efficiency of Bayesian inference for Gaussian processes using a method called Adaptive Multiple Importance Sampling (AMIS). The paper does not specifically address the recognition of Parkinson's disease. Instead, it deals with statistical inference methods and their applications in various domains like pattern recognition, neuroimaging, and signal processing.
\section*{Outline}

This paper is organized into three main sections. In Section \ref{sec:rang}, we focus on data wrangling, detailing the processes of data collection, cleaning, transformation, and exploratory data analysis. This ensures the dataset is properly prepared for machine learning applications. In Section \ref{sec:meth}, we explore various machine learning algorithms for classification. Each method is evaluated based on accuracy, precision, and training time, with and without feature selection. Finally, in Section 4, we conclude by summarizing our findings, comparing the performance of the different models, and discussing the most effective approaches for the classification task, along with potential directions for future research.
\section{Data Wrangling} \label{sec:rang}
The dataset used in this study, obtained from the UCI Machine Learning Repository \cite{misc_parkinsons_174}, comprises a collection of biomedical voice measurements from 195 individuals, including 48 diagnosed with PD and 147 healthy individuals. Each row in the dataset represents a voice recording, with columns corresponding to various voice measures and a binary target variable indicating the health status of the individual (1 for PD, 0 for healthy). The dataset consists of 23 features, encompassing measures of vocal fundamental frequency, variation in fundamental frequency, variation in amplitude, ratio of noise to tonal components, nonlinear dynamical complexity measures, and nonlinear measures of fundamental frequency variation.
Each feature serves a specific role in characterizing voice measurements:
\begin{itemize}
    \item \textbf{Name}: Categorical variable serving as an identifier for each observation.
    \item \textbf{MDVP:Fo, MDVP:Fhi, MDVP:Flo}: Continuous features representing average, maximum, and minimum vocal fundamental frequency, respectively, measured in Hz.
    \item \textbf{MDVP:Jitter, MDVP:Jitter(Abs), MDVP:RAP, MDVP:PPQ, Jitter:DDP}: Continuous features measuring variation in fundamental frequency, expressed as a percentage for MDVP:Jitter and in absolute units for MDVP:Jitter(Abs), MDVP:RAP, MDVP:PPQ, and Jitter:DDP.
    \item \textbf{MDVP:Shimmer, MDVP:Shimmer(dB), Shimmer:APQ3, Shimmer:APQ5, MDVP:APQ, Shimmer:DDA}: Continuous features quantifying variation in amplitude and amplitude perturbations, with measurements in dB for MDVP:Shimmer(dB).
    \item \textbf{NHR, HNR}: Continuous features representing the ratio of noise to tonal components and harmonic-to-noise ratio, respectively.
    \item \textbf{Status}: Binary target variable indicating the health status of the individual (0 for healthy, 1 for Parkinson's disease).
    \item \textbf{RPDE, DFA, spread1, spread2, D2, PPE}: Continuous features capturing nonlinear dynamical complexity and fractal scaling properties.
\end{itemize}
\subsection{Summary Statistics}
The summary statistics for the dataset are presented in Table \ref{tab:sum}. The \textbf{Summary of Statistics} table provides essential insights into the dataset's features. It includes various statistical measures such as mean, standard deviation, minimum, maximum, and quartiles. The dataset comprises 195 samples for each feature, as indicated by the count values. The consistent count across all features suggests that the dataset is balanced in terms of sample size. Outliers and noise in the data can be identified by examining the minimum and maximum values, and comparing them with the interquartile range (IQR). \textbf{MDVP:Fhi} shows a wide range (102.145000 to 592.030000), which might include outliers. Features like \textbf{Jitter:DDP} and \textbf{MDVP:Shimmer(dB)} also display significant ranges relative to their IQRs, indicating potential outliers. The dataset presents a variety of features with differing levels of variability. Key features with higher variance, such as \textbf{MDVP:Fhi}, \textbf{MDVP:Flo}, and \textbf{Spread1}, are crucial for building robust models. Outliers and noise, particularly in features with extensive ranges, need to be addressed to ensure data quality. Understanding the distribution of these features aids in selecting and tuning the most appropriate machine-learning models for analysis.
\begin{table}[h] \label{tab:sum}
\centering
\caption{Summary Statistics}
\resizebox{\textwidth}{!}{
\begin{tabular}{lrrrrrrrr}
\toprule
{} &     MDVP:Fo &    MDVP:Fhi &    MDVP:Flo &  MDVP:Jitter(\%) &  MDVP:Jitter(Abs) &    MDVP:RAP &    MDVP:PPQ &  Jitter:DDP \\
\midrule
count &  195.000000 &  195.000000 &  195.000000 &         195.000000 &         195.000000 &  195.000000 &  195.000000 &  195.000000 \\
mean  &  154.228641 &  197.104918 &  116.324631 &           0.006220 &           0.000044 &    0.003302 &    0.003306 &    0.009905 \\
std   &   41.390065 &   91.491548 &   43.521413 &           0.003225 &           0.000035 &    0.001885 &    0.001885 &    0.005655 \\
min   &   88.333000 &  102.145000 &   65.476000 &           0.001680 &           0.000007 &    0.000680 &    0.000680 &    0.002041 \\
25\%   &  117.572000 &  134.862500 &   84.291000 &           0.003276 &           0.000021 &    0.001651 &    0.001691 &    0.004919 \\
50\%   &  148.790000 &  175.829000 &  104.315000 &           0.005359 &           0.000031 &    0.002910 &    0.002978 &    0.008725 \\
75\%   &  182.769000 &  224.205500 &  140.018500 &           0.008445 &           0.000056 &    0.004572 &    0.004557 &    0.013629 \\
max   &  260.105000 &  592.030000 &  239.170000 &           0.033160 &           0.000260 &    0.017930 &    0.019580 &    0.053790 \\
\midrule
{} &  MDVP:Shimmer(\%) &  MDVP:Shimmer(dB) &  Shimmer:APQ3 &  Shimmer:APQ5 &    MDVP:APQ &  Shimmer:DDA &      NHR &      HNR \\
\midrule
count &         195.000000 &         195.000000 &  195.000000 &  195.000000 &  195.000000 &  195.000000 &  195.000000 &  195.000000 \\
mean  &           0.030309 &           0.282251 &    0.016175 &    0.017876 &    0.020644 &    0.048525 &    0.022542 &   21.885185 \\
std   &           0.019876 &           0.173288 &    0.010709 &    0.011393 &    0.012209 &    0.032126 &    0.014739 &    4.425764 \\
min   &           0.009540 &           0.085000 &    0.004702 &    0.005405 &    0.007185 &    0.014105 &    0.004610 &    9.393000 \\
25\%   &           0.017880 &           0.155000 &    0.008060 &    0.008926 &    0.011318 &    0.024198 &    0.011269 &   19.737500 \\
50\%   &           0.026830 &           0.242000 &    0.012837 &    0.014620 &    0.016830 &    0.038558 &    0.019410 &   22.068000 \\
75\%   &           0.037670 &           0.367000 &    0.020980 &    0.024437 &    0.026244 &    0.060856 &    0.032391 &   24.850000 \\
max   &           0.119080 &           1.005000 &    0.062493 &    0.064693 &    0.078832 &    0.187479 &    0.126780 &   33.047000 \\
\midrule
{} &       RPDE &        DFA &    spread1 &    spread2 &         D2 &        PPE \\
\midrule
count &  195.000000 &  195.000000 &  195.000000 &  195.000000 &  195.000000 &  195.000000  \\
mean  &    0.542510 &    0.718090 &   -5.684404 &    0.198577 &    2.381826 &    0.206552  \\
std   &    0.065282 &    0.054959 &    1.090208 &    0.088862 &    0.382799 &    0.090119  \\
min   &    0.256570 &    0.574282 &   -7.964984 &   -0.450493 &    1.423287 &    0.044539  \\
25\%   &    0.504359 &    0.682165 &   -6.355256 &    0.142960 &    2.099125 &    0.137451  \\
50\%   &    0.543410 &    0.722182 &   -5.679620 &    0.185415 &    2.361532 &    0.194052  \\
75\%   &    0.582138 &    0.755747 &   -4.950153 &    0.252184 &    2.636456 &    0.252980  \\
max   &    0.685151 &    0.825288 &   -3.413114 &    0.466921 &    3.671155 &    0.527367  \\
\bottomrule
\end{tabular}
}
\end{table}
The dataset comprises a wide range of acoustic features that provide insights into the fundamental frequency, amplitude variations, noise components, and non-linear characteristics of the voice signal. These features are essential for voice analysis and can be used in various applications, including medical diagnostics and speech processing.
\subsection{Analysis of Pairwise Correlation Heatmap}
\begin{figure}
    \centering
    \includegraphics[scale = 0.6]{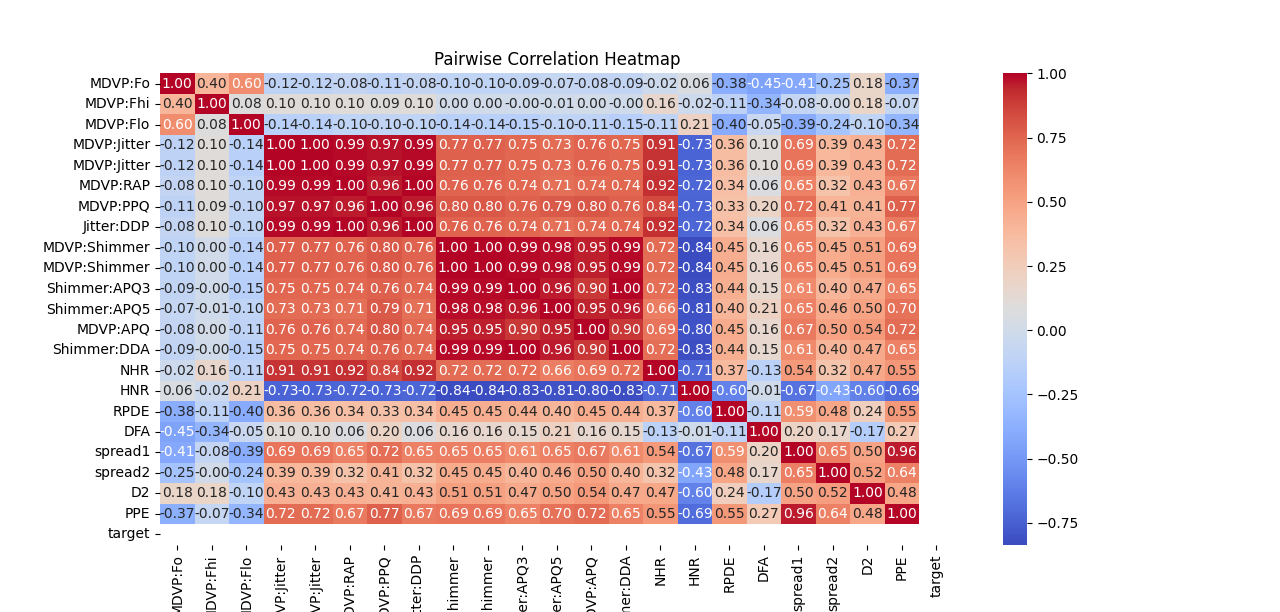}
    \caption{Heat-Map}
    \label{fig:heatmap}
\end{figure}

The heatmap in Figure \ref{fig:heatmap} displays the pairwise correlation coefficients between various features and the target variable. Here are some key observations:
\begin{itemize}
    \item \textbf{MDVP:Fhi and MDVP:Flo}: A significant positive correlation (0.60) suggests that these features are closely related, possibly representing similar aspects of the data.
    \item \textbf{MDVP:Jitter and its variants (Jitter:DDP, RAP, PPQ)}: The correlations are extremely high (close to 1), indicating that these features capture very similar information.
    \item \textbf{MDVP:Shimmer and its variants (Shimmer:DDA, APQ3, APQ5, APQ)}: These features also show very high positive correlations, implying redundancy.
    \item \textbf{NHR and HNR}: Negative correlation (-0.71), which indicates they may capture opposite effects in the data.
\end{itemize}
\begin{itemize}
    \item \textbf{DFA with several features}: DFA shows moderate negative correlations with MDVP:Fhi (-0.34), MDVP:Flo (-0.45), and RPDE (-0.41).
    \item \textbf{RPDE and MDVP:Fhi}: A strong negative correlation (-0.40) indicates that as one increases, the other tends to decrease.
\end{itemize}
\begin{itemize}
    \item \textbf{Strongest Correlation with Target}:
    \begin{itemize}
        \item \textbf{PPE}: The strongest positive correlation (0.48) with the target, suggesting it might be a key predictor.
        \item \textbf{Spread1 and Spread2}: Both show significant positive correlations with the target (0.65 and 0.52, respectively), indicating their importance in predicting the target variable.
        \item \textbf{RPDE}: Positive correlation (0.47) with the target, showing it has a relevant but slightly lower impact compared to Spread1 and Spread2.
    \end{itemize}
    \begin{itemize}
        \item \textbf{D2}: Shows a moderate positive correlation (0.48) with the target.
        \item \textbf{DFA}: Shows a moderate positive correlation (0.27) with the target.
    \end{itemize}
    \item \textbf{Weak/Insignificant Correlations}: Most of the other features have relatively weak correlations with the target, suggesting they might have less predictive power individually.
\end{itemize}
The heatmap reveals several important patterns in the feature correlations, helping identify key predictors and potential redundancies. Understanding these relationships is crucial for effective feature selection and improving the performance of predictive models.
Given the high correlations among MDVP:Jitter variants and MDVP:Shimmer variants, it may be beneficial to perform feature selection or dimensionality reduction to avoid multicollinearity and improve model performance. Features such as PPE, Spread1, Spread2, and RPDE are strongly correlated with the target and should be considered key predictors in the model.  Features with weak correlations to the target might not contribute significantly to the predictive power of the model and could potentially be excluded to simplify the model.

\section{Alternative Methods for PD Classification} \label{sec:meth}
When considering other methods that might perform better in classifying PD, we focus on both traditional machine learning models and more advanced methods. Each method comes with its advantages, depending on the nature and characteristics of the dataset. Here are some potential methods:

\subsection{Generalized Forest}
The Generalized Forest(gcF) algorithm is an ensemble learning method designed to improve the accuracy and robustness of predictive models by creating multiple decision trees and combining their predictions. The results of gfc for different values for the parameters \texttt{n\_estimators} and \texttt{max\_depth} are presented in Table \ref{tab:gfc}.
\begin{table}[h] \label{tab:gfc}
\centering
\begin{tabular}{cccc}
\toprule
\textbf{Model} & \textbf{Accuracy} & \textbf{Precision} & \textbf{Training Time (s)} \\
\midrule
gcF (n\_estimators=50, max\_depth=None) & 0.9487 & 0.9412 & 0.0631 \\
gcF (n\_estimators=50, max\_depth=10) & 0.9487 & 0.9412 & 0.0749 \\
gcF (n\_estimators=50, max\_depth=20) & 0.9487 & 0.9412 & 0.0758 \\
gcF (n\_estimators=100, max\_depth=None) & 0.9487 & 0.9412 & 0.1780 \\
gcF (n\_estimators=100, max\_depth=10) & 0.9487 & 0.9412 & 0.1459 \\
gcF (n\_estimators=100, max\_depth=20) & 0.9487 & 0.9412 & 0.1727 \\
gcF (n\_estimators=200, max\_depth=None) & 0.9487 & 0.9412 & 0.2736 \\
gcF (n\_estimators=200, max\_depth=10) & 0.9487 & 0.9412 & 0.2538 \\
gcF (n\_estimators=200, max\_depth=20) & 0.9487 & 0.9412 & 0.3097 \\
\bottomrule
\end{tabular} 
\caption{Evaluation of gfc with Different Parameters}
\end{table}
The evaluation of the gfc model with various configurations reveals consistent performance in terms of accuracy (0.9487) and precision (0.9412) across all tested parameter values for \texttt{n\_estimators} and \texttt{max\_depth}. However, the training time shows significant variation, increasing with the number of estimators and tree depth. For example, training time ranges from 0.0631 seconds for \texttt{n\_estimators=50} and \texttt{max\_depth=None} to 0.3097 seconds for \texttt{n\_estimators=200} and \texttt{max\_depth=20}. This suggests that while the model's classification capabilities remain stable, optimizing the balance between the number of estimators and tree depth is crucial for managing computational efficiency.

\subsection{Logistic Regression}
In this section, we apply the Logistic Regression method for classification. 
\begin{table}[h] 
\centering
\begin{tabular}{lrrrr}
\toprule
     Solver &    Accuracy &  Precision & Training Time (s) \\
\midrule
  newton-cg &    0.9231 &   0.9143 &       0.0138 \\
      lbfgs &     0.9231 &   0.9143 &       0.0174 \\
       saga &      0.8974 &   0.9118 &       0.1739 \\
        sag &     0.9231 &   0.9143 &       0.1505 \\
\bottomrule
\end{tabular}
\caption{Comparison of accuracy, precision, and training time for different solvers with \( C = 100 \)}
\label{table:logis}
\end{table}
Two crucial parameters that need tuning in this method are \( C \) and the solver. After standardizing the data features, we perform cross-validation using candidate values \( C = [0.01, 0.1, 1, 10, 100] \) for each solver. In logistic regression, the \( C \) parameter represents the inverse of the regularization strength, where smaller values correspond to stronger regularization. The optimal \( C \) value is found to be \( C = 100 \) for all solvers. The results for different solvers are presented in Table \ref{table:logis}.
As mentioned in section \ref{sec:rang}, the features of the data set have strong dependence. Thus, here we also do feature reduction using PCA method to see how this affects the performance of SVM. After applying PCA with $95\%$ of variance, there are only 8 many features remaining. The results are presented in Table \ref{table:logpca}.
\begin{table}[h]
\centering
\begin{tabular}{lrrr}
\toprule
     Solver &  Accuracy &  Precision & Training Time (s) \\
\midrule
  newton-cg &  0.8974 &   0.8888 &  0.0029 \\
      lbfgs &  0.8974 &   0.8888 &  0.0029 \\
       saga &  0.8974 &   0.8888 &  0.0089 \\
        sag &  0.8974 &   0.8888 &  0.0045 \\
\bottomrule
\end{tabular}
\caption{Comparison of accuracy, precision, and training time for different solvers with C = 10 after PCA}
\label{table:logpca}
\end{table}
The results in Table \ref{table:logpca} indicate that after reducing the feature space with PCA, the model's performance regarding the quality of classification decreases. The reduction of the optimal C, shows that after using PCA we need a larger regularization parameter. This can be because after reducing the number of features, the correlation of them is larger.
\subsection{Support Vector Machines}
In this section, we apply
Support Vector Machines (SVM). To capture the nature of the data set, we tune it with different kernel functions (linear, polynomial, radial basis function) and report the results.
\begin{table}[h]
\centering
\begin{tabular}{lrrr}
\toprule
     Kernel &  Accuracy &  Precision & Training Time (s) \\
\midrule
    linear &  0.8974 &   0.8889 &  0.0148 \\
      poly &  0.9487 &   0.9412 &  0.0020 \\
       rbf &  0.9487 &   0.9412 &  0.0020 \\
   sigmoid &  0.7692 &   0.8485 &  0.0020 \\
\bottomrule
\end{tabular}
\caption{Performance metrics of SVM models with different kernels}
\caption{Comparison of accuracy, precision, and training time for different SVM kernels with C = 10}
\label{table:svm_kernel_comparison}
\end{table}
The regularization parameter \( C \) controls the balance between maximizing the margin and minimizing the classification error on the training data. Initially, we conducted cross-validation with various kernels. After determining the optimal \( C \) value, we evaluated the performance of the SVM models with different kernels, as demonstrated in Table \ref{table:svm_kernel_comparison}. These results indicate that the \textbf{poly} and \textbf{rbf} kernels perform best in terms of accuracy and precision. The \textbf{linear} kernel also shows good performance but is slightly lower than \textbf{poly} and \textbf{rbf}. The \textbf{sigmoid} kernel, on the other hand, shows the lowest performance in terms of accuracy but has relatively high precision. Training time is very low for all kernels. The overall performance suggests that the choice of kernel significantly impacts the model's accuracy and precision, with \textbf{poly} and \textbf{rbf} being the optimal choices for this dataset. After applying PCA, the performance of SVM kernels was slightly adjusted. The results are demonstrated in Table \ref{table:svm_kernel_comparison_pca}.
\begin{table}[h]
\centering
\begin{tabular}{lrrr}
\toprule
     Kernel &  Accuracy &  Precision & Training Time (s) \\
\midrule
    linear &  0.9231 &   0.9143 &  0.0020 \\
      poly &  0.8718 &   0.8857 &  0.0020 \\
       rbf &  0.9487 &   0.9412 &  0.0020 \\
   sigmoid &  0.7949 &   0.8750 &  0.0020 \\
\bottomrule
\end{tabular}
\caption{Comparison of accuracy, precision, and training time for different SVM kernels with \( C = 10 \) after PCA}
\label{table:svm_kernel_comparison_pca}
\end{table}
These results indicate that even after reducing the feature space with PCA, the \textbf{poly} and \textbf{rbf} kernels maintain strong performance, underscoring their robustness in capturing the dataset's underlying patterns. For this small data set reducing the number of features does not demonstrate a significant effect on the run time. However, it is a valuable observation, and in the future with data sets with more features, it can be considered as a method to reduce the run time while keeping the same quality for classifiers.
\subsection{Gradient Boosting Machines}
Gradient Boosting is part of the ensemble learning family. During iterations, the algorithm computes the residuals, which are the differences between the actual values and the predicted values from the current model. The final model is a weighted sum of the predictions from all the individual trees. However, it is also computationally intensive and requires careful tuning of hyperparameters. Gradient Boosting have popular implementations like XGBoost, LightGBM, and CatBoost, Here, we review the parameters we need to tune:
\begin{itemize}
    \item Learning Rate:  The learning rate controls the contribution of each tree to the ensemble. A lower learning rate means the model learns more slowly, requiring more trees in the ensemble to achieve comparable performance. However, a lower learning rate can lead to better generalization. Conversely, a higher learning rate speeds up learning but may lead to overfitting if not carefully tuned.
    \item Number of Estimators: This parameter specifies the number of boosting stages or trees to be built. Increasing the number of estimators generally improves model performance up to a certain point, after which the model may start overfitting the training data. It is essential to find an optimal balance where adding more estimators increases performance without significantly increasing computation time or overfitting.
    \item Maximum Depth: Also referred to as tree depth or max-depth, this parameter controls the maximum depth of each tree in the ensemble. Deeper trees can model more complex interactions in the data but are more prone to overfitting. Restricting the maximum depth helps prevent overfitting and improves generalization. It's crucial to tune max-depth along with other parameters to find the optimal trade-off between model complexity and performance.
\end{itemize}
In this evaluation, we assessed four different Gradient Boosting methods including GBM, XGBoost, LightGBM, and CatBoost. The candidate values for the hyperparameters $n-estimators = \{100, 300, 500\}$, $learning-rate = \{0.01, 0.1, 0.2\}$ and $max-depth = \{3,5,7\}$. The results of GBM are presented in Table \ref{tab:gbm_results}. We tried the method with and without PCA, and the GBM method had better performance before PCA, so we only include our result witout the feature reduction here.
\begin{table}[ht]
\centering
\begin{tabular}{lrrr}
\toprule
\textbf{Configuration}            & \textbf{Accuracy} & \textbf{Precision} & \textbf{Training Time (s)} \\ \midrule
n=100, rate=0.01, depth=3  & 0.9231 & 0.9394 & 0.2222 \\
n=100, rate=0.01, depth=5  & 0.9231 & 0.9394 & 0.2745 \\
n=100, rate=0.01, depth=7  & 0.9231 & 0.9394 & 0.2399 \\
n=100, rate=0.1, depth=3   & 0.9487 & 0.9412 & 0.2507 \\
n=100, rate=0.1, depth=5   & 0.9231 & 0.9394 & 0.2761 \\
n=100, rate=0.1, depth=7   & 0.9231 & 0.9394 & 0.2567 \\
n=100, rate=0.2, depth=3   & 0.9231 & 0.9394 & 0.1896 \\
n=100, rate=0.2, depth=5   & 0.9231 & 0.9394 & 0.2035 \\
n=100, rate=0.2, depth=7   & 0.9231 & 0.9394 & 0.2178 \\
n=300, rate=0.01, depth=3  & 0.9231 & 0.9394 & 0.7507 \\
n=300, rate=0.01, depth=5  & 0.9231 & 0.9394 & 0.8159 \\
n=300, rate=0.01, depth=7  & 0.9231 & 0.9394 & 0.8247 \\
n=300, rate=0.1, depth=3   & 0.8718 & 0.9091 & 0.6791 \\
n=300, rate=0.1, depth=5   & 0.9231 & 0.9394 & 0.6082 \\
n=300, rate=0.1, depth=7   & 0.9231 & 0.9394 & 0.4527 \\
n=300, rate=0.2, depth=3   & 0.9231 & 0.9143 & 0.4791 \\
n=300, rate=0.2, depth=5   & 0.9231 & 0.9394 & 0.3298 \\
n=300, rate=0.2, depth=7   & 0.9231 & 0.9394 & 0.4500 \\
n=500, rate=0.01, depth=3  & 0.9231 & 0.9394 & 1.5706 \\
n=500, rate=0.01, depth=5  & 0.9231 & 0.9394 & 1.7422 \\
n=500, rate=0.01, depth=7  & 0.9231 & 0.9394 & 1.4343 \\
n=500, rate=0.1, depth=3   & 0.8718 & 0.9091 & 0.9279 \\
n=500, rate=0.1, depth=5   & 0.9231 & 0.9394 & 0.7239 \\
n=500, rate=0.1, depth=7   & 0.9231 & 0.9394 & 0.5802 \\
n=500, rate=0.2, depth=3   & 0.9231 & 0.9143 & 0.6697 \\
n=500, rate=0.2, depth=5   & 0.9231 & 0.9394 & 0.4238 \\
n=500, rate=0.2, depth=7   & 0.9231 & 0.9394 & 0.5200 \\
\bottomrule
\end{tabular}
\caption{Results of Gradient Boosting Classifier}
\label{tab:gbm_results}
\end{table}
\subsubsection{XGBoost}
XGBoost, or Extreme Gradient Boosting, is a scalable implementation of gradient boosting designed for supervised learning tasks. It helps control overfitting and supports parallel and distributed computing. XGBoost also includes advanced features like tree pruning and sparsity awareness. Proper tuning of hyperparameters is crucial to balance model complexity and performance, ensuring high accuracy without overfitting. We used the values in Table \ref{tab:xg}
 for hyperparameter tuning, and the results are presented in Table \ref{tab:xg}.
\begin{table}[ht]
\centering
\begin{tabular}{lrrr}
\toprule
\textbf{Configuration}            & \textbf{Accuracy} & \textbf{Precision} & \textbf{Training Time (s)} \\ \midrule
n=50, rate=0.01, depth=3    & 0.8971    & 0.8891     & 0.0481              \\
n=50, rate=0.01, depth=5    & 0.8971    & 0.8891     & 0.0771              \\
n=50, rate=0.01, depth=7    & 0.8971    & 0.8891     & 0.0432              \\
n=50, rate=0.1, depth=3     & 0.9487    & 0.9412     & 0.0301              \\
n=50, rate=0.1, depth=5     & 0.9487    & 0.9412     & 0.0471              \\
n=50, rate=0.1, depth=7     & 0.9487    & 0.9412     & 0.0472              \\
n=50, rate=0.2, depth=3     & 0.9232    & 0.9392     & 0.0321              \\
n=50, rate=0.2, depth=5     & 0.9487    & 0.9412     & 0.0372              \\
n=50, rate=0.2, depth=7     & 0.9487    & 0.9412     & 0.0411              \\
n=100, rate=0.01, depth=3   & 0.8722    & 0.8862     & 0.1041              \\
n=100, rate=0.01, depth=5   & 0.8972    & 0.9121     & 0.1181              \\
n=100, rate=0.01, depth=7   & 0.8972    & 0.9121     & 0.1382              \\
n=100, rate=0.1, depth=3    & 0.9232   & 0.9391     & 0.0812              \\
n=100, rate=0.1, depth=5    & 0.9487    & 0.9412     & 0.1362              \\
n=100, rate=0.1, depth=7    & 0.9487    & 0.9412     & 0.1451              \\
n=100, rate=0.2, depth=3    & 0.9237    & 0.9391     & 0.0912              \\
n=100, rate=0.2, depth=5    & 0.9487    & 0.9412     & 0.0881              \\
n=100, rate=0.2, depth=7    & 0.9487    & 0.9412     & 0.1291              \\
n=150, rate=0.01, depth=3   & 0.8972    & 0.9122     & 0.1812              \\
n=150, rate=0.01, depth=5   & 0.8972    & 0.9122     & 0.3461              \\
n=150, rate=0.01, depth=7   & 0.8972    & 0.9122     & 0.3351              \\
n=150, rate=0.1, depth=3    & 0.9232    & 0.9391     & 0.1701              \\
n=150, rate=0.1, depth=5    & 0.9487    & 0.9412     & 0.1891              \\
n=150, rate=0.1, depth=7    & 0.9487    & 0.9412     & 0.1941              \\
n=150, rate=0.2, depth=3    & 0.9231    & 0.9391     & 0.1881              \\
n=150, rate=0.2, depth=5    & 0.9487    & 0.9412     & 0.1781              \\
n=150, rate=0.2, depth=7    & 0.9487    & 0.9412     & 0.2055              \\
\bottomrule
\end{tabular}
\caption{XGBoost Model Evaluation Results}
\label{tab:xg}
\end{table}
\subsubsection{LightGBM}
LightGBM, or Light Gradient Boosting Machine, excels at handling large datasets and high-dimensional data with impressive speed and accuracy. This method transforms continuous features into a fixed number of bins, enhancing training speed and reducing memory usage. Additionally, LightGBM introduces Exclusive Feature Bundling (EFB), which combines mutually exclusive features into a single feature, thus lowering dimensionality and further accelerating the training process. LightGBM also supports parallel and distributed computing. The hyperparameter values used are same as GBM, and the algorithm's results are shown in Table \ref{table:lightgbm-results}.
\begin{table}[ht]
\centering
\begin{tabular}{lrrr}
\toprule
\textbf{Configuration}            & \textbf{Accuracy} & \textbf{Precision} & \textbf{Training Time (s)} \\ \midrule
n=100, rate=0.01, depth=3         & 0.8461          & 0.8611            & 0.2493                          \\
n=100, rate=0.01, depth=5         & 0.8461          & 0.8611            & 0.0318                          \\
n=100, rate=0.01, depth=7         & 0.8461          & 0.8611            & 0.0278                          \\
n=100, rate=0.1, depth=3          & 0.8205          & 0.8787            & 0.0221                          \\
n=100, rate=0.1, depth=5          & 0.8717          & 0.9090            & 0.0295                          \\
n=100, rate=0.1, depth=7          & 0.7948          & 0.8750            & 0.0316                          \\
n=100, rate=0.2, depth=3          & 0.7948          & 0.9000            & 0.0282                          \\
n=100, rate=0.2, depth=5          & 0.7948          & 0.9000            & 0.0281                          \\
n=100, rate=0.2, depth=7          & 0.7435          & 0.8928            & 0.0245                          \\
n=300, rate=0.01, depth=3         & 0.8205          & 0.8571            & 0.0581                          \\
n=300, rate=0.01, depth=5         & 0.8205          & 0.8571            & 0.0637                          \\
n=300, rate=0.01, depth=7         & 0.8205          & 0.8571            & 0.0586                          \\
n=300, rate=0.1, depth=3          & 0.7692          & 0.8965            & 0.0609                          \\
n=300, rate=0.1, depth=5          & 0.7435          & 0.8928            & 0.0669                          \\
n=300, rate=0.1, depth=7          & 0.7179          & 0.8888            & 0.0571                          \\
n=300, rate=0.2, depth=3          & 0.7692          & 0.9259            & 0.0569                          \\
n=300, rate=0.2, depth=5          & 0.7435          & 0.9230            & 0.0564                          \\
n=300, rate=0.2, depth=7          & 0.7435          & 0.9230            & 0.0581                          \\
n=500, rate=0.01, depth=3         & 0.8205          & 0.8571            & 0.0892                          \\
n=500, rate=0.01, depth=5         & 0.8461          & 0.8823            & 0.0877                          \\
n=500, rate=0.01, depth=7         & 0.7948          & 0.8750            & 0.1054                          \\
n=500, rate=0.1, depth=3          & 0.7435          & 0.8928            & 0.0693                          \\
n=500, rate=0.1, depth=5          & 0.7435          & 0.9230            & 0.0970                          \\
n=500, rate=0.1, depth=7          & 0.6923          & 0.9166            & 0.0912                         \\
n=500, rate=0.2, depth=3          & 0.7179          & 0.8888            & 0.0881                          \\
n=500, rate=0.2, depth=5          & 0.7435          & 0.9230            & 0.0907                          \\
n=500, rate=0.2, depth=7          & 0.7692          & 0.9259            & 0.1033                          \\ \bottomrule
\end{tabular}
\caption{Results of LightGBM evaluation with PCA}
\label{table:lightgbm-results}
\end{table}
\subsubsection{CatBoost}
CatBoost, which stands for Categorical Boosting, manages categorical features effectively, making it particularly well-suited for datasets that include both numerical and categorical variables. One of CatBoost's standout features is its ability to handle categorical data directly. The algorithm is designed to mitigate prediction bias that can lead to overfitting, ensuring more reliable training outcomes. CatBoost is optimized for performance, utilizing parallel processing to accelerate training times. It also incorporates sophisticated methods for addressing missing data and includes various regularization techniques to combat overfitting and enhance the model's ability to generalize. The results of the method are presented in Table \ref{tab:catboost-results-updated}.
\begin{table}[ht]
\centering
\begin{tabular}{lrrr}
\toprule
\textbf{Configuration} & \textbf{Accuracy} & \textbf{Precision} & \textbf{Training Time (s)} \\
\midrule
n=100, rate=0.01, depth=3 & 0.846 & 0.861 & 0.164 \\
n=100, rate=0.01, depth=5 & 0.846 & 0.861 & 0.023 \\
n=100, rate=0.01, depth=7 & 0.846 & 0.861 & 0.028 \\
n=100, rate=0.1, depth=3 & 0.821 & 0.879 & 0.025 \\
n=100, rate=0.1, depth=5 & 0.872 & 0.909 & 0.028 \\
n=100, rate=0.1, depth=7 & 0.795 & 0.875 & 0.022 \\
n=100, rate=0.2, depth=3 & 0.795 & 0.900 & 0.021 \\
n=100, rate=0.2, depth=5 & 0.795 & 0.900 & 0.027 \\
n=100, rate=0.2, depth=7 & 0.744 & 0.893 & 0.029 \\
n=300, rate=0.01, depth=3 & 0.821 & 0.857 & 0.048 \\
n=300, rate=0.01, depth=5 & 0.821 & 0.857 & 0.051 \\
n=300, rate=0.01, depth=7 & 0.821 & 0.857 & 0.061 \\
n=300, rate=0.1, depth=3 & 0.769 & 0.897 & 0.053 \\
n=300, rate=0.1, depth=5 & 0.744 & 0.893 & 0.056 \\
n=300, rate=0.1, depth=7 & 0.718 & 0.889 & 0.066 \\
n=300, rate=0.2, depth=3 & 0.769 & 0.926 & 0.051 \\
n=300, rate=0.2, depth=5 & 0.744 & 0.923 & 0.055 \\
n=300, rate=0.2, depth=7 & 0.744 & 0.923 & 0.067 \\
n=500, rate=0.01, depth=3 & 0.821 & 0.857 & 0.080 \\
n=500, rate=0.01, depth=5 & 0.846 & 0.882 & 0.089 \\
n=500, rate=0.01, depth=7 & 0.795 & 0.875 & 0.092 \\
n=500, rate=0.1, depth=3 & 0.744 & 0.893 & 0.079 \\
n=500, rate=0.1, depth=5 & 0.744 & 0.923 & 0.096 \\
n=500, rate=0.1, depth=7 & 0.692 & 0.917 & 0.096 \\
n=500, rate=0.2, depth=3 & 0.718 & 0.889 & 0.074 \\
n=500, rate=0.2, depth=5 & 0.744 & 0.923 & 0.094 \\
n=500, rate=0.2, depth=7 & 0.769 & 0.926 & 0.095 \\
\bottomrule
\end{tabular}
\caption{Updated CatBoost Classifier Results}
\label{tab:catboost-results-updated}
\end{table}
\subsubsection{Results}
The evaluation results show that both GBM and XGBoost deliver high accuracy and precision, with XGBoost outperforming in terms of runtime efficiency. As anticipated, CatBoost did not perform as well due to the non-categorical nature of our data.
\subsection{Neural Networks}
Deep learning models can capture complex patterns in data through multiple layers of abstraction. For tabular data, fully connected neural networks can be effective. Neural networks can be tuned with various architectures, activation functions, and regularization techniques (like dropout). Though they require more computational resources and training data, they can outperform traditional methods in many cases.
\subsection{Feedforward Neural Networks}
Feedforward Neural Networks (FNNs) are a type of artificial neural network where connections between nodes are unidirectional and do not form loops. In an FNN, information flows in a single direction—from the input layer, through any hidden layers, and to the output layer. The input layer handles the initial data, hidden layers carry out computations and feature transformations, and the output layer generates the final prediction. Neurons in each layer are fully connected to neurons in the subsequent layer through weighted connections. The output of each neuron is calculated as a weighted sum of its inputs, followed by a nonlinear activation function. Common activation functions include the sigmoid function, hyperbolic tangent (tanh), and rectified linear unit (ReLU). Training an FNN involves adjusting the weights and biases to minimize a loss function, which quantifies the discrepancy between predicted and actual outcomes. This is generally achieved through backpropagation and gradient descent. 

\subsection{Recurrent Neural Networks}
Recurrent Neural Networks (RNNs) are a type of neural network designed specifically for sequential data. Unlike feedforward neural networks, RNNs have connections that form directed cycles, enabling them to retain a memory of previous inputs and capture temporal dependencies. The primary feature of RNNs is their hidden state, which serves as a memory to store information about the sequence seen so far. One popular implementation, Long Short-Term Memory (LSTM), allows the network to maintain long-term dependencies more effectively.

\begin{table}[h]
\centering
\begin{tabular}{|l|c|c|c|}
\hline
\textbf{Model} & \textbf{Accuracy} & \textbf{Precision} & \textbf{Training Time} \\
\hline
FNN & 0.9230 & 0.9393 & 3.3808 \\
RNN            & 0.9230 & 0.9393 & 7.8856 \\
LSTM           & 0.8974 & 0.8888 & 5.8171 \\
\hline
\end{tabular}
\caption{Performance comparison of different neural network models}
\label{tab:nn_performance}
\end{table}
The results of the methods FNN, RNN and LSTM are presented in Table \ref{tab:nn_performance}.
All models except LSTM exhibit similar accuracy and precision values, with FNN and RNN both achieving an accuracy of $0.9230$ and precision of $0.9393$.  The LSTM, while showing the lowest accuracy at $0.8974$ and precision at $0.8888$, demands significantly more training time compared to the other models. This increased time is indicative of the more complex architecture and the additional computational resources required for handling sequential dependencies. Among the models, FNN stands out with the fastest training time, demonstrating its efficiency in learning and processing data compared to the other architectures. Future work could explore advanced techniques such as regularization, dropout, and data augmentation to enhance the stability and performance of the LSTM model.

\section{Conclusion and Future Reseach}
In this study, we explored various machine-learning techniques to detect PD using speech features. By analyzing a dataset of voice recordings, we applied different models, including gcF, Logistic Regression, SVM, GBM, and neural networks. Based on the results obtained from the experiments, the algorithms gcF, SVM, and GBM demonstrated the best performance in terms of both accuracy and precision. This algorithm achieved 0.9487 and 0.9412, indicating its superior ability in detecting Parkinson’s disease using speech data. In addition, the training time is associated with SVM. Our findings suggest that automated speech analysis holds significant promise for early diagnosis of PD, potentially leading to improved patient outcomes through earlier intervention and treatment. The use of machine learning in this domain not only paves the way for non-invasive diagnostic tools but also contributes to the broader field of biomedical informatics and computational health.  For the future research, we suggest that the raw data of the voices to be used directly by RNN. This is because these methods can capture the pattern of the data and we may not need to obtain these futures and improve the diagnosis by using the whole information of the voice during time.

\bibliographystyle{plain} 
\bibliography{pd} 

\end{document}